\documentclass[francais]{gretsi}

\usepackage[utf8]{inputenc}

\usepackage[T1]{fontenc}

\usepackage{amsmath, amssymb, amsfonts}
\usepackage{amsthm}
\usepackage{stmaryrd}

\usepackage[font=footnotesize,skip=3pt]{caption}
\usepackage[font=small,skip=0pt]{subcaption}

\theoremstyle{definition}

\usepackage{setspace}

\usepackage{xcolor}

\newcommand{\smallparagraph}[1]{\medskip\noindent\textbf{#1}\hspace{1em}}

\newcommand\matseq[2]{\mathbf{#1}^{(#2)}}
\newcommand\bflysupp[1]{\mathbf{S}_{\mathrm{bf}}^{(#1)}}
\newcommand\set[1]{\llbracket#1\rrbracket}
\newcommand\supp{\mathrm{supp}}
\newcommand\identity[1]{\mathbf{I}_{#1}}
\newcommand\monarch[1]{\mathcal{M}^{(#1)}}
\newcommand\rowpartition[1]{R^{(#1)}}
\newcommand\colpartition[1]{C^{(#1)}}
\newcommand{\bflytx}{T^X_{\mathrm{bf}}}
\newcommand{\bflyto}{T^\Omega_{\textrm{bf}}}
\newcommand\bflyerror[3]{E_{\textrm{bf}}(#1, #2, #3)}

\newcommand\extrafootertext[1]{%
    \bgroup
    \renewcommand\thefootnote{\fnsymbol{footnote}}%
    \renewcommand\thempfootnote{\fnsymbol{mpfootnote}}%
    \footnotetext[0]{#1}%
    \egroup
}

\begin{document}

\titre{Factorisation butterfly par identification algorithmique\\de blocs de rang un}

\auteurs{
  \auteur{Léon}{Zheng}{leon.zheng@valeo.com}{1,2}
  \auteur{Gilles}{Puy}{gilles.puy@valeo.com}{2}
  \auteur{Elisa}{Riccietti}{elisa.riccietti@ens-lyon.fr}{1}
  \auteur{Patrick}{Pérez}{patrick.perez@valeo.com}{2}
  \auteur{Rémi}{Gribonval}{remi.gribonval@inria.fr}{1}
}

\affils{
  \affil{1}{Univ Lyon, EnsL, UCBL, CNRS, Inria, LIP, F-69342, LYON Cedex 07, France.
  }
  \affil{2}{valeo.ai, Paris, France.
  }
}

\resume{
Plusieurs matrices associées à des transformées rapides possèdent une certaine propriété de rang faible qui se caractérise par l'existence de plusieurs partitions par blocs de la matrice, où chaque bloc est de rang faible. \`A condition de connaître ces partitions, il existe alors des algorithmes, dits de \emph{factorisation butterfly}, qui approchent la matrice en un produit de facteurs creux, permettant ainsi une évaluation rapide de l'opérateur linéaire associé.
Cet article propose une nouvelle méthode pour identifier algorithmiquement les partitions en blocs de rang faible d'une matrice admettant une factorisation butterfly, sans hypothèse analytique sur ses coefficients. 

}

\abstract{
Many matrices associated with fast transforms posess a certain low-rank property characterized by the existence of several block partitionings of the matrix, where each block is of low rank. Provided that these partitionings are known, there exist algorithms, called \emph{butterfly factorization} algorithms, that approximate the matrix into a product of sparse factors, thus enabling a rapid evaluation of the associated linear operator. 
This paper proposes a new method to identify algebraically these block partitionings for a matrix admitting a butterfly factorization, without any analytical assumption on its entries. 
}

\maketitle

\section{Introduction}
L'évaluation rapide d'un opérateur linéaire est un enjeu clé dans de nombreux domaines comme le calcul scientifique, le traitement du signal ou l'apprentissage automatique. 
Dans des applications mettant en jeu un très %
grand nombre de paramètres, le calcul \emph{direct} de la multiplication matrice-vecteur passe difficilement à l'échelle pour cause de %
complexité quadratique en la taille de la matrice. 
De nombreux travaux se sont ainsi intéressés à la construction d'algorithmes rapides pour la multiplication matricielle, en s'appuyant typiquement sur des propriétés analytiques ou algébriques des matrices apparaissant dans les problèmes étudiés.

\begin{figure}[h!]
    \centering
     \begin{subfigure}[b]{0.23\linewidth}
         \centering
         \includegraphics[width=\linewidth]{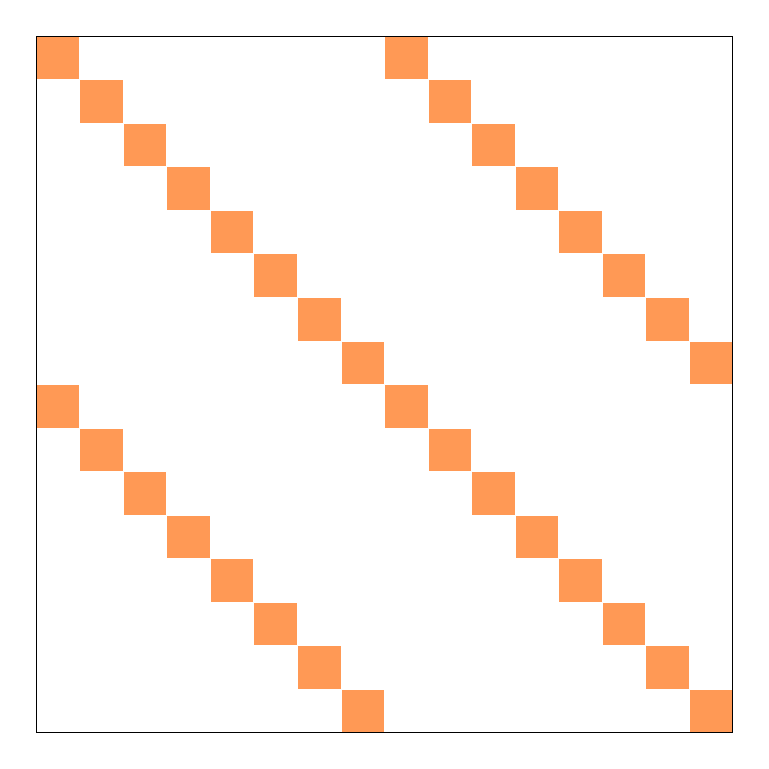}
         \caption{$\bflysupp{1}$}
     \end{subfigure}
     \hfill
     \begin{subfigure}[b]{0.23\linewidth}
         \centering
         \includegraphics[width=\linewidth]{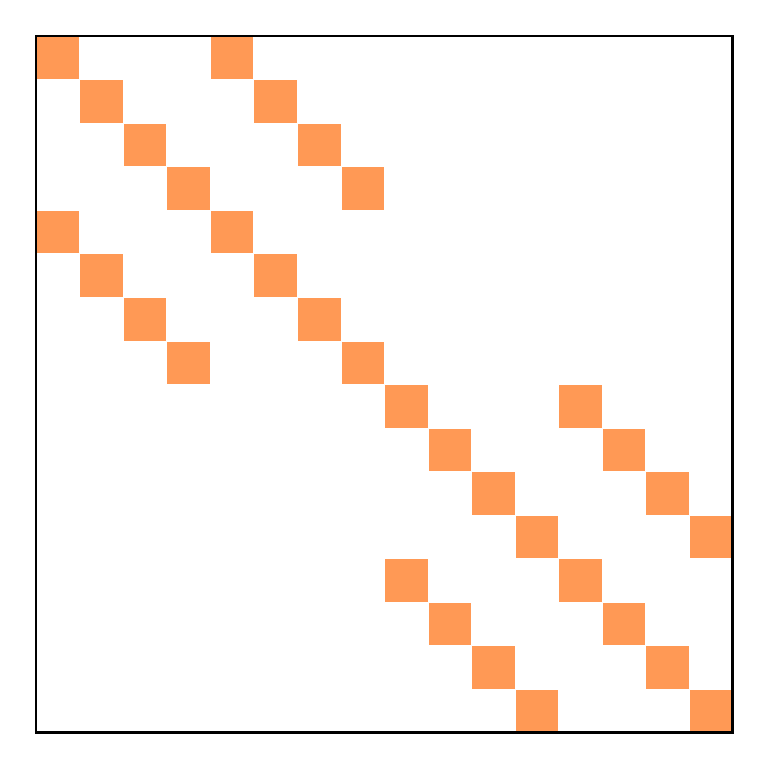}
         \caption{$\bflysupp{2}$}
     \end{subfigure}
          \hfill
      \begin{subfigure}[b]{0.23\linewidth}
         \centering
         \includegraphics[width=\linewidth]{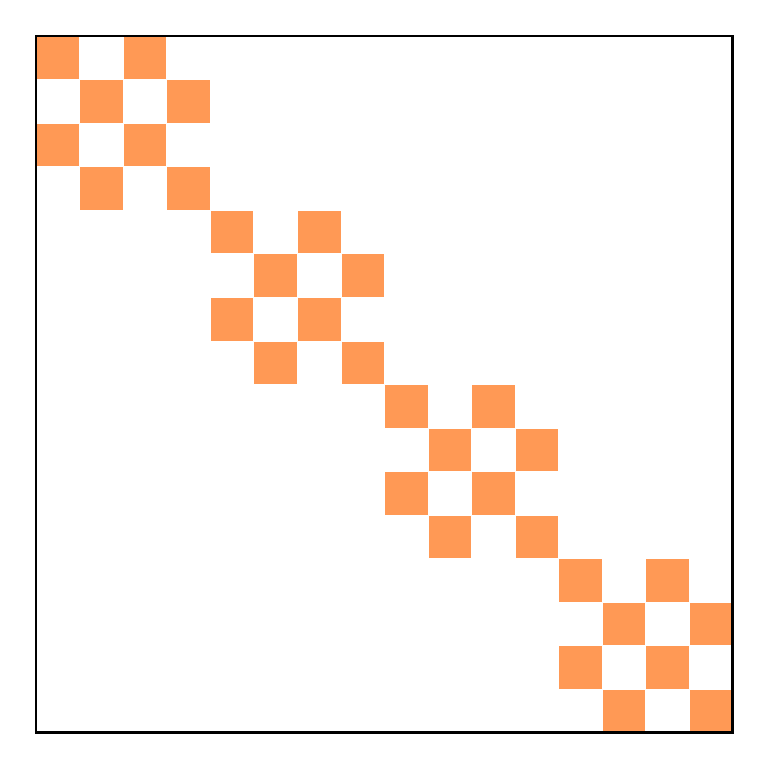}
         \caption{$\bflysupp{3}$}
     \end{subfigure}
          \hfill
      \begin{subfigure}[b]{0.23\linewidth}
         \centering
         \includegraphics[width=\linewidth]{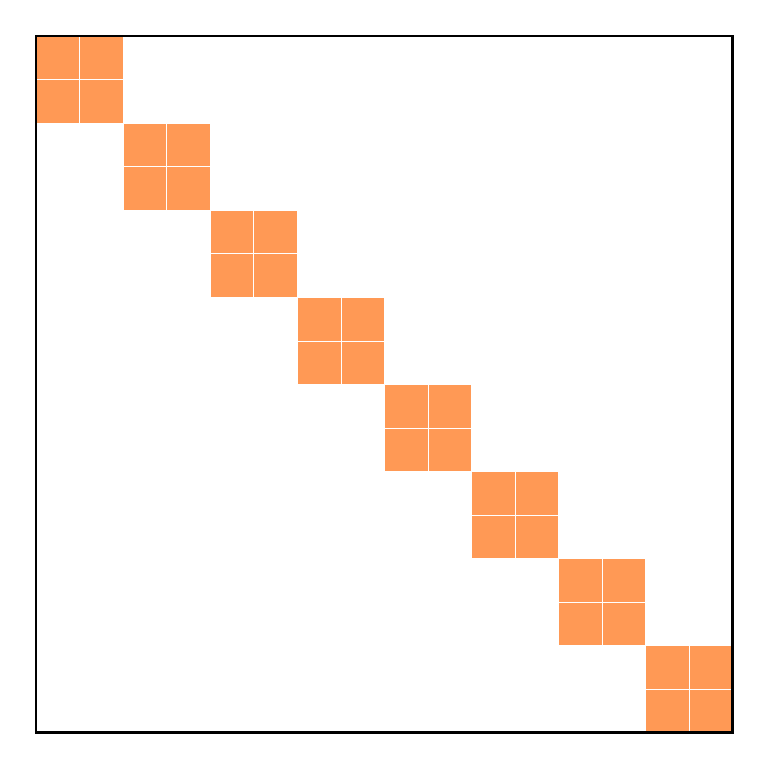}
         \caption{$\bflysupp{4}$}
     \end{subfigure}
     \caption{Supports des facteurs butterfly, $N=16$ (cf.~section \ref{sec:formulation}).}
    \label{fig:support-butterfly}
\end{figure}

Ces algorithmes rapides sont souvent associés à une factorisation creuse de la matrice correspondante, comme c'est le cas pour la matrice de Hadamard ou  de la transformée de Fourier discrète. En effet, à permutation près des lignes et des colonnes, ces matrices de taille $N$ possèdent une factorisation \emph{butterfly}, dans le sens où elle s'écrivent comme le produit de $\mathcal{O}(\log N)$ facteurs ayant chacun $\mathcal{O}(N)$ coefficients non nuls, et dont les supports ont une structure particulière illustrée dans la figure~\ref{fig:support-butterfly}. Il a été montré dans \cite{dao2019learning} que la classe des matrices admettant une telle factorisation est expressive, au sens où elle %
contient plusieurs matrices structurées utilisées en apprentissage ou en traitement du signal. 
Ce modèle serait donc pertinent pour chercher des factorisations creuses d'opérateurs pour lesquels un algorithme 
d'évaluation rapide n'est pas connu.

Trouver un algorithme rapide associé à la factorisation butterfly se formalise alors comme un problème d'optimisation, où l'on minimise en norme de Frobenius $\| \cdot \|_F$ l’erreur d’approximation d'une matrice $\mathbf{A} \in \mathbb{C}^{N\times N}$ par un produit de facteurs butterfly $\matseq{X}{1}, \ldots, \matseq{X}{L}$, à permutations près des lignes et colonnes, encodées par les matrices de permutation $\mathbf{P}, \mathbf{Q}$:
\begin{equation}
\label{eq:bfly-fact}
    \min_{(\matseq{X}{\ell})_{\ell=1}^L, \mathbf{P}, \mathbf{Q}} \| \mathbf{A} - \mathbf{Q}^\top \matseq{X}{1} \ldots  \matseq{X}{L} \mathbf{P} \|_F.
\end{equation}

\begin{figure}[t!]
    \centering
    \includegraphics[width=\linewidth]{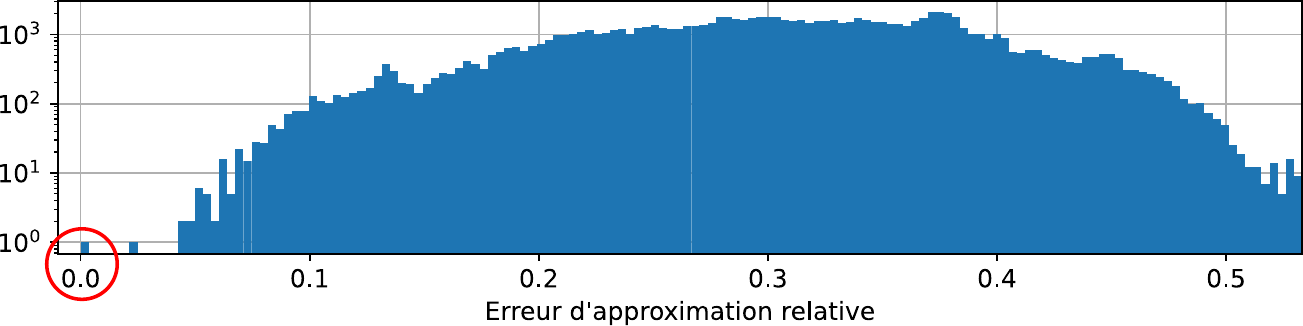}
    \caption{Histogramme des erreurs relatives pour la résolution de \eqref{eq:bfly-fact} dans le cas sans bruit et $N=8$, en énumérant et en fixant chaque permutation de lignes et de colonnes, comme expliqué en section \ref{sec:nécessité}. 
    Parmi toutes les paires de permutation de lignes et colonnes, seule une paire donne une erreur nulle, comme indiqué par le cercle rouge.
    }
    \label{fig:recherche-exhaustive}
\end{figure}

Lorsque les permutations optimales $\mathbf{P}, \mathbf{Q}$ sont connues, %
il existe un algorithme hiérarchique efficace en complexité $\mathcal{O}(N^2)$  pour trouver des facteurs butterfly $(\matseq{X}{\ell})_{\ell=1}^L$ donnant une faible erreur d'approximation, avec garanties de reconstruction dans le cas du problème sans bruit \cite{le2022fast,zheng2023efficient}. 
Mais lorsque ces permutations ne sont pas connues, le problème \eqref{eq:bfly-fact} est conjecturé comme étant difficile : d'une part, si l'on énumère toutes les permutations possibles pour résoudre \eqref{eq:bfly-fact}, on constate numériquement que seule une petite proportion de permutations donne une faible erreur d'approximation, comme illustré en figure \ref{fig:recherche-exhaustive}, ce qui montre la \emph{nécessité} d'identifier les bonnes permutations afin de résoudre \eqref{eq:bfly-fact} ; mais d'autre part, une recherche exhaustive de toutes les permutations n'est pas tractable, même en tenant compte de certaines équivalences de permutations vis-à-vis de \eqref{eq:bfly-fact}, comme discuté en section \ref{sec:formulation}.

Afin d'identifier les permutations optimales, nous nous appuyons sur le fait qu'une matrice admettant une factorisation butterfly possède une certaine propriété dite de \emph{rang faible complémentaire} \cite{li2015butterfly}, dans le sens où il existe des partitions par blocs de la matrice, où chaque bloc est de rang faible. 
Il suffit alors d'identifier ces partitions pour résoudre le problème \eqref{eq:bfly-fact}. 

Ceci peut se faire \emph{analytiquement} lorsque les coefficients de la matrice s'expriment par un noyau régulier $(\mathbf{x}, \boldsymbol{\omega}) \mapsto K(\mathbf{x}, \boldsymbol{\omega})$ évalué sur des paramètres $\{ \mathbf{x}_i \}_{i=1}^N$, $\{ \boldsymbol{\omega}_j \}_{j=1}^N$, par exemple pour %
des matrices %
associées à certains %
opérateurs intégraux \cite{candes2009fast} ou transformées spéciales de fonctions \cite{o2010algorithm}.

En revanche, si la matrice étudiée n'a pas de forme analytique, ou si %
celle-ci n'est pas accessible, la littérature ne propose pas, à notre connaissance, de méthode pour identifier ces partitions. Aussi, %
nous proposons une heuristique à base de partitionnement spectral alterné des lignes et des colonnes pour identifier les partitions en blocs de rang faible, sans hypothèse analytique. En termes d'applications, cette heuristique permet de vérifier algorithmiquement qu'un opérateur linéaire, pour lequel on ignore l'existence d'un algorithme d'évaluation rapide, possède une propriété de rang faible complémentaire approchée %
donnant lieu à une bonne approximation de la matrice associée par un produit de facteurs butterfly.

La section \ref{sec:formulation} détaille cette propriété de rang faible utilisée pour identifier les permutations optimales. La section \ref{sec:nécessité} montre numériquement la nécessité d'identifier ces permutations, ce qui motive notre méthode expliquée en section \ref{sec:méthode}, et validée empiriquement en section \ref{sec:expérience}.

\section{Formulation du problème}
\label{sec:formulation}
Dans le reste du papier, nous nous plaçons dans le cas des matrices carrées de taille $N := 2^L$ pour un certain entier $L \geq 2$. Pour une matrice quelconque $\mathbf{M}$, son support $\supp(\mathbf{M})$ est l'ensemble des indices correspondants aux coefficients non nuls de $\mathbf{M}$. Notons $\identity{N}$ la matrice identité de taille $N$, et $\otimes$ le produit de Kronecker. Selon \cite{dao2019learning,le2022fast,zheng2023efficient}, $\mathbf{A} \in \mathbb{C}^{N \times N}$ est une \emph{matrice butterfly} si elle admet une factorisation $\mathbf{A} = \matseq{X}{1} \ldots \matseq{X}{L}$, où chaque facteur $\matseq{X}{\ell} \in \mathbb{C}^{N \times N}$, nommé \emph{facteur butterfly}, satisfait la contrainte de support fixe $\supp(\matseq{X}{\ell}) \subseteq \supp(\bflysupp{\ell})$ pour $\ell \in \set{L} := \{1, \ldots, L\}$, avec $\bflysupp{\ell} := \identity{2^{\ell-1}} \otimes [\begin{smallmatrix}
  1 & 1\\
  1 & 1
\end{smallmatrix} ] \otimes \identity{N / 2^\ell}$. L'ensemble des matrices admettant une telle factorisation est noté $\mathcal{B}$. Les supports $\bflysupp{\ell}$ sont illustrés par la figure \ref{fig:support-butterfly} et, par abus de notation, la contrainte de support fixe sera notée $\supp(\matseq{X}{\ell}) \subseteq \bflysupp{\ell}$ par la suite. 

Le problème que nous souhaitons résoudre est donc \eqref{eq:bfly-fact}, sous la contrainte $\supp(\matseq{X}{\ell}) \subseteq \bflysupp{\ell}$ pour tout $\ell \in \set{L}$.
Pour trouver une bonne solution de l'instance du problème où les matrices de permutations $\mathbf{P}$ et $\mathbf{Q}$ sont fixées, on applique l'algorithme hiérarchique de \cite{le2022fast,zheng2023efficient}. Par la suite, on note $\bflyerror{\mathbf{A}}{\mathbf{P}}{\mathbf{Q}} := \| \mathbf{A} - \mathbf{Q}^\top \tilde{\mathbf{X}}^{(1)} \ldots \tilde{\mathbf{X}}^{(L)} \mathbf{P} \|_F$ l'erreur d'approximation donnée par la sortie $(\tilde{\mathbf{X}}^{(1)}, \ldots, \tilde{\mathbf{X}}^{(L)})$ de l'algorithme \ref{algo:hierarchical}, quand $\mathbf{P}$ et $\mathbf{Q}$ sont fixées. En revanche, le cas difficile est celui où les permutations ne sont pas fixées : le reste de la section explique 
notre approche qui s'appuie sur la propriété de rang faible complémentaire.

\LinesNumbered
\begin{algorithm}
    \Entree{$\mathbf{A} \in \mathbb{C}^{N \times N}$, permutations $\mathbf{P}, \mathbf{Q}$.}
  Résoudre $\min \| \mathbf{Q} \mathbf{A} \mathbf{P}^\top - \matseq{X}{1} \ldots  \matseq{X}{L} \|_F$ sous la contrainte %
  $\supp(\matseq{X}{\ell}) \subseteq \bflysupp{\ell}$ pour tout $\ell \in \set{L}$ via \cite[Algorithme 3.2]{zheng2023efficient} \\
  \Retour{les facteurs butterfly $(\tilde{\mathbf{X}}^{(1)}, \ldots, \tilde{\mathbf{X}}^{(L)})$}
  \caption{Algorithme hiérarchique pour résoudre l'instance de \eqref{eq:bfly-fact} où $\mathbf{P}$ et $\mathbf{Q}$ sont fixées.}
  \label{algo:hierarchical}
\end{algorithm}

\smallparagraph{Propriété de rang faible complémentaire \cite{li2015butterfly}}
Celle-ci se définit à l'aide de deux \og cluster tree \fg $T^X$ et $T^\Omega$, qui sont des arbres binaires avec $L=\log_2(N)$ niveaux (sans compter la racine) dont chaque n{\oe}ud est un sous-ensemble non vide de $\set{N}$, avec comme racine $\set{N}$ au niveau 0, et où les enfants constituent une partition de leur parent en deux sous-ensembles de même cardinal \cite{hackbusch2015hierarchical}.
Une matrice $\mathbf{A}$ possède la propriété de rang faible complémentaire pour $T^X$ et $T^\Omega$ si, pour chaque $\ell \in \set{L-1}$, pour chaque n{\oe}ud $R$ au niveau $L - \ell$ de $T^X$ et pour chaque n{\oe}ud $C$ au niveau de $\ell$ de $T^\Omega$, la restriction $\mathbf{A}_{R, C}$ de $\mathbf{A}$ sur les lignes et colonnes indexées par $R$ et $C$ est de rang faible. 
Sachant que les n{\oe}uds d'un même niveau de $T^X$ et ceux de $T^\Omega$ forment respectivement une partition des indices de lignes et de colonnes,
la propriété de rang faible complémentaire impose que les blocs $\{ \mathbf{A}_{R, C} \}_{R, C}$ des partitions décrites par les niveaux de $T^X$ et $T^\Omega$ soient de rang faible.

Il a été montré dans \cite{zheng2023efficient} que toute matrice  butterfly $\mathbf{A} \in \mathcal{B}$ satisfait la propriété de rang faible complémentaire (de rang 1) pour une paire d'arbres $(\bflytx, \bflyto)$ spécifique, comme expliqué ci-dessous pour rappeler les concepts qui nous seront utiles par la suite.
Notons $\bflysupp{p:q} := \bflysupp{p} \ldots \bflysupp{q}$ ($1 \leq p < q \leq L$), et définissons la $\ell$-ème classe Monarch \cite{dao2022monarch} pour $\ell \in \set{L-1}$ comme étant l'ensemble
\begin{equation*}
    \monarch{\ell} := \left \{ \mathbf{X} \mathbf{Y}, \, \supp(\mathbf{X}) \subseteq \bflysupp{1:\ell}, \, \supp(\mathbf{Y}) \subseteq \bflysupp{\ell+1:L} \right \}.
\end{equation*}
On vérifie alors que $\mathcal{B} \subseteq \bigcap_{\ell=1}^{L-1} \monarch{\ell}$, et qu'une matrice $\mathbf{A}$ appartient à $\monarch{\ell}$ si, et seulement si, $\mathbf{A}_{R, C}$ est de rang au plus 1 pour tout $(R, C) \in \mathcal{P}^{(\ell)} := \{ \rowpartition{\ell}_i \}_{i=1}^{N/2^\ell} \times \{ \colpartition{\ell}_j \}_{j=1}^{2^\ell}$, où
\begin{align*}
    \rowpartition{\ell}_i &:= \{ i + (k-1) N / 2^\ell, \, k \in \set{2^\ell} \}, \\
    \colpartition{\ell}_j &:= \{ (j-1) N / 2^\ell + k, \, k \in \set{N/2^\ell} \}.
\end{align*}
On observe que $\{ \rowpartition{\ell}_i \}_{i=1}^{N/2^\ell}$ et $\{ \colpartition{\ell}_j \}_{j=1}^{2^\ell}$ constituent chacun une partition de $\set{N}$ : ainsi, lorsque $\mathbf{A} \in \monarch{\ell}$, ces partitions de lignes et de colonnes décrivent une partition de $\mathbf{A}$ en blocs de rang 1.
On définit alors $\bflytx$ et $\bflyto$ comme étant les deux arbres pour lesquels les n{\oe}uds de $\bflytx$ au niveau $L - \ell$ et les n{\oe}uds de $\bflyto$ au niveau $\ell$ sont précisément $\{ \rowpartition{\ell}_i \}_{i=1}^{N/2^\ell}$ et $\{ \colpartition{\ell}_j \}_{j=1}^{2^\ell}$, pour chaque $\ell \in \set{L-1}$. 
Ainsi, $\mathbf{A} \in \bigcap_{\ell=1}^{L-1} \monarch{\ell}$ est précisément une reformulation de la propriété de rang faible complémentaire (de rang 1) pour les arbres $\bflytx$ et $\bflyto$.

\smallparagraph{Permutation d'indices}
Pour tout \og cluster tree \fg $T$ dont la racine est $\set{N}$, et pour toute permutation $\sigma: \set{N} \to \set{N}$, on définit le \og cluster tree \fg $\sigma(T)$ obtenu en permutant selon $\sigma$ les indices dans $T$. En particulier, on remarque que pour tout arbre $T^X$ et $T^\Omega$, il existe plusieurs matrices de permutation $\mathbf{P}$ et $\mathbf{Q}$ pour lesquelles $\sigma_{\mathbf{Q}}(\bflytx) = T^X$ et $\sigma_{\mathbf{P}}(\bflyto) = T^\Omega$, où $\sigma_{\mathbf{P}}, \sigma_{\mathbf{Q}}$ sont les permutations associées aux matrices $\mathbf{P}$, $\mathbf{Q}$. Ceci définit alors des classes d'équivalence de permutations de lignes $[\mathbf{Q}_{T^X}]$ et de colonnes $[\mathbf{P}_{T^\Omega}]$. 

\smallparagraph{Approche pour résoudre \eqref{eq:bfly-fact}} Supposons que la matrice cible dans \eqref{eq:bfly-fact} soit de la forme $\mathbf{A} := \tilde{\mathbf{Q}}^\top \tilde{\mathbf{A}} \tilde{\mathbf{P}}$, où $\tilde{\mathbf{A}} \in \mathcal{B}$ est une matrice butterfly, et $\tilde{\mathbf{P}}$, $\tilde{\mathbf{Q}}$ sont deux matrices de permutations arbitraires inconnues. Puisque $\tilde{\mathbf{A}}$ satisfait la propriété de rang faible complémentaire pour les arbres $\bflytx$ et $\bflyto$, la matrice $\tilde{\mathbf{Q}}^\top \tilde{\mathbf{A}} \tilde{\mathbf{P}}$ la satisfait également mais pour les arbres $T^X := \sigma_{\tilde{\mathbf{Q}}}(\bflytx)$ et $T^\Omega := \sigma_{\tilde{\mathbf{P}}}(\bflyto)$.  Si l'on parvient à reconstruire $T^X$ et $T^\Omega$ à partir de l'observation de $\mathbf{A}$, alors \eqref{eq:bfly-fact} peut se résoudre en choisissant une paire quelconque $(\mathbf{P}, \mathbf{Q}) \in [\mathbf{P}_{T^\Omega}] \times [\mathbf{Q}_{T^X}]$, et en appliquant l'algorithme \ref{algo:hierarchical} avec ces permutations fixées. En effet, un tel choix suffit pour garantir que $\mathbf{Q} \mathbf{A} \mathbf{P}^\top$ satisfasse la propriété de rang faible complémentaire pour $\bflytx$ et $\bflyto$, i.e., $\mathbf{Q} \mathbf{A} \mathbf{P}^\top \in \bigcap_{\ell=1}^{L-1} \monarch{\ell}$, puisque 
les matrices $\mathbf{P}^\top$ et $\mathbf{Q}^\top$ sont associées aux permutations inverses
$\sigma^{-1}_{\mathbf{P}}$ et $\sigma_{\mathbf{Q}}^{-1}$, et par définition des classes d'équivalence, on a bien $\sigma_{\mathbf{P}}^{-1}(T^\Omega) = \bflyto$ et $\sigma^{-1}_{\mathbf{Q}}(T^X) = \bflytx$.
En conclusion, afin de résoudre \eqref{eq:bfly-fact}, il est suffisant d'identifier les arbres $T^X$ et $T^\Omega$ pour lesquelles la matrice cible satisfait la propriété de rang faible complémentaire, ce qui se ramène, par définition, à identifier des partitions de $\mathbf{A}$ en blocs de rang 1.

\section{Nécessité de retrouver les partitions}
\label{sec:nécessité}

Nous montrons à présent empiriquement  que  l'identification des arbres $T^X$ et $T^\Omega$ est en fait nécessaire pour résoudre \eqref{eq:bfly-fact} avec $\mathbf{A} := \tilde{\mathbf{Q}}^\top \tilde{\mathbf{A}} \tilde{\mathbf{P}}$. 
Considérons $\tilde{\mathbf{A}} \in \mathcal{B}$ dont les facteurs ont des coefficients non nuls tirés selon une gaussienne centrée réduite. Puis, nous énumérons tous les arbres $T^X, T^\Omega$ possibles, calculons une solution via l'algorithme \ref{algo:hierarchical} en fixant une paire arbitraire $(\mathbf{P}, \mathbf{Q}) \in [\mathbf{P}_{T^\Omega}] \times [\mathbf{Q}_{T^X}]$,
et vérifions que les seuls arbres donnant une erreur $\bflyerror{\mathbf{A}}{\mathbf{P}}{\mathbf{Q}}$ faible sont $T^X$ et $T^\Omega$.
Un dénombrement de tous les arbres montre en revanche que cette expérience n'est pas tractable pour une grande taille $N$. En effet, le nombre $u_N$ de \og cluster tree \fg pour une même racine de cardinal $N$ satisfait la relation de récurrence $u_N = \frac{1}{2}
    {N \choose \frac{N}{2}}
    (u_{\frac{N}{2}})^2$ avec $u_2 = 1$, car il y a $\frac{1}{2}
    {N \choose \frac{N}{2}}$ paires d'enfants possibles pour la racine, et chaque enfant est un \og cluster tree \fg dont la racine est de cardinal $N/2$.
Ainsi, nous considérons $N=8$ à titre illustratif, ce qui donne $u_8 = 315$.

Dans la figure \ref{fig:recherche-exhaustive}, 
la recherche exhaustive sur toutes les paires d'arbres montrent que l'erreur est nulle seulement pour \emph{une seule} paire, et que les autres paires échouent à la résolution de \eqref{eq:bfly-fact}. Ceci illustre donc empiriquement la nécessité d'identifier les bons arbres $T^X$ et $T^\Omega$ pour résoudre \eqref{eq:bfly-fact}. Démontrer formellement une telle nécessité pour n'importe quelle taille de matrices pourra faire l'objet de futurs travaux.

\section{Partitionnement spectral alterné}
\label{sec:méthode}

Nous proposons
l'algorithme \ref{algo:heuristic} à base de partitionnement spectral alterné pour identifier les arbres $T^X$ et $T^\Omega$ pour lesquels la matrice $\mathbf{A}$ satisfait la propriété de rang faible complémentaire. 
Pour décrire notre approche, %
commençons par expliquer la résolution du problème $\min_{\mathbf{M} \in \monarch{\ell}, \mathbf{P}, \mathbf{Q}} \| \mathbf{A} - \mathbf{Q}^\top \mathbf{M} \mathbf{P} \|_F^2$ pour chaque $\ell \in \set{L-1}$.
\'Etant donné qu'une matrice $\mathbf{M}$ appartient à $\monarch{\ell}$ si et seulement si le bloc $\mathbf{M}_{R, C}$ est de rang au plus 1 pour tout jeu de lignes et colonnes $(R, C) \in \mathcal{P}^{(\ell)}$, ce problème équivaut à
\begin{equation}
    \label{eq:min-partition}
    \min_{\{ R_i \}_{i=1}^{N/2^\ell}, \{ C_j \}_{j=1}^{2^\ell}} \sum_{i=1}^{N/2^\ell} \sum_{j=1}^{2^\ell} \min_{\mathbf{x}, \mathbf{y}} \|\mathbf{A}_{R_i, C_j} - \mathbf{x} \mathbf{y}^*\|_F^2,
\end{equation}
où $\{ R_i \}_{i=1}^{N/2^\ell}$, $\{ C_j \}_{j=1}^{2^\ell}$ sont des partitions de lignes et de colonnes en sous-ensembles de même cardinal, et $\min_{\mathbf{x}, \mathbf{y}} \|\mathbf{A}_{R, C} - \mathbf{x} \mathbf{y}^*\|_F^2$ calcule la meilleure approximation de rang 1 de $\mathbf{A}_{R, C}$. Le symbole $*$ désigne la matrice adjointe.

\LinesNumbered
\begin{algorithm}[t]
\Entree{$\mathbf{A}, \ell, \alpha > 0$, graine aléatoire $\texttt{seed}$.}
    $\{ C_j \}_{j=1}^{2^\ell} \gets$ partition aléatoire de $\set{N}$ selon $\texttt{seed}$ \\
\Pour{$i=1,\ldots, n$}{
  $\{ R_i \}_{i=1}^{N/2^\ell} \gets$ partitionnement en fixant $\{ C_j \}_{j=1}^{\ell}$ \\
    $\{ C_j \}_{j=1}^{2^\ell} \gets$ partitionnement en fixant $\{ R_i \}_{i=1}^{N/2^\ell}$
}
    $E_\ell \gets \sum_{i,j} \min_{\mathbf{x}, \mathbf{y}} \|\mathbf{A}_{R_i, C_j} - \mathbf{x} \mathbf{y}^*\|_F^2$ \\
    \Sortie{$E_\ell, \{ R_i \}_{i=1}^{N/2^\ell}, \{ C_j \}_{j=1}^{2^\ell}$.}
  \caption{Partitionnement spectral alterné pour résoudre \eqref{eq:min-partition} durant $n$ itérations, à $\ell \in \set{L-1}$ fixé.}
    \label{algo:partitionnement}
\end{algorithm}

\smallparagraph{Partitionnement des lignes~}
Pour qu'une optimisation alternée fonctionne, il est nécessaire de pouvoir résoudre le problème \ref{eq:min-partition} lorsqu'une des deux partitions est connue. 
Fixons ainsi sans perte de généralité une partition de colonnes $\{ C_j \}_{j=1}^{2^\ell}$, et cherchons une partition de lignes $\{ R_i \}_{i=1}^{N/2^\ell}$ qui minimise \eqref{eq:min-partition}. 

Pour cela, nous nous inspirons des méthodes existantes pour le problème de <<\,subspace clustering\,>> \cite{vidal2016generalized}. Définissons $2^\ell$ graphes $\{ \mathcal{G}_j \}_{j=1}^{2^\ell}$, où les $N$ n{\oe}uds du graphe $\mathcal{G}_j$ sont les $N$ lignes de $\mathbf{A}$ restreintes aux colonnes $C_j$ --notées $\mathbf{A}_{k, C_j}$ pour $k \in \set{N}$-- et les poids des arêtes de $\mathcal{G}_j$ sont donnés par la matrice de similarité $\mathbf{W}^{(j)} \in \mathbb{R}^{N \times N}$ définie par
\begin{equation*}
    \mathbf{W}^{(j)}_{k, l} := \left( \frac{| {\mathbf{A}_{k, C_j}}^* \mathbf{A}_{l, C_j} |}{\| \mathbf{A}_{k, C_j} \|_2 \|\mathbf{A}_{l, C_j} \|_2} \right)^\alpha \quad \forall k, l \in \set{N},
\end{equation*}
avec $\alpha > 0$ un paramètre contrôlant le contraste entre les poids des arêtes. 
Intuitivement, un groupe de lignes restreintes aux colonnes $C_j$ est inter-connecté par des poids de fortes valeurs dans le graphe $\mathcal{G}_j$ lorsque les lignes correspondantes sont corrélées, ce qui est le cas lorsqu'elles forment un bloc de rang 1. Inversement, deux lignes non corrélées ont un poids de faible valeur.
Ainsi, en résolvant un problème de coupe minimale sur le graphe $\mathcal{G}$ dont la matrice de similarité est $\mathbf{W} := \sum_{j=1}^{2^\ell} \mathbf{W}^{(j)}$, les groupes de n{\oe}uds obtenus doivent correspondre à un partitionnement des lignes qui minimise \eqref{eq:min-partition}.

Concrètement, à partir de la matrice de similarité $\mathbf{W}$, un partitionnement spectral \cite{von2007tutorial} du graphe $\mathcal{G}$ est effectué en calculant la décomposition en vecteurs propres du Laplacien non normalisé $\mathbf{L} := \mathbf{D} - \mathbf{W}$, où $\mathbf{D}$ est la matrice de degrés dont les coefficients diagonaux sont 
$\mathbf{W} (1 \ldots 1
)^\top$. Afin de garantir un partitionnement des lignes en groupe de même taille, l'étape de partitionnement $k$-moyenne sur les représentations spectrales est implémentée selon la méthode de \cite{bradley2000constrained}.

\smallparagraph{Optimisation alternée~}
Quand la partition de colonnes n'est plus fixée, la résolution de \eqref{eq:min-partition} suit l'algorithme \ref{algo:partitionnement} de partitionnement \emph{alterné}, où l'on initialise aléatoirement une partition des colonnes et, à chaque itération, un partitionnement spectral des lignes est effectuée en fixant la partition de colonnes de l'itération précédente, et vice-versa en échangeant le rôle des lignes et des colonnes. %
Sans garanties de réussite, l'algorithme \ref{algo:partitionnement} peut nécessiter plusieurs réinitialisations pour trouver une 
solution.

\LinesNumbered
\begin{algorithm}[t]
\Entree{$\mathbf{A}, \{\alpha_k \}_{k=1}^K, \{ \texttt{seed}_m \}_{m=1}^M$.}
\Pour{$\ell=1,\ldots, L-1$}{
    \Pour{$k = 1, \ldots, K$, $m=1, \ldots, M$}
    {
        $E_\ell, \{ \tilde{R}_i^{(\ell)} \}_i, \{  \tilde{C}_j^{(\ell)} \}_j \gets $ algorithme \ref{algo:partitionnement} appliqué à $(\mathbf{A}, \ell, \alpha_k, \texttt{seed}_m)$
    }
    Garder $\{ \tilde{R}_i^{(\ell)} \}_i, \{ \tilde{C}_j^{(\ell)} \}_j$ avec le plus petit $E_\ell$ 
}
    Vérifier que $\{ \tilde{R}_i^{(\ell)} \}_i$, $\{ \tilde{C}_j^{(\ell)} \}_j$ pour $\ell \in \set{L-1}$ forment chacun des \og cluster tree \fg $T^X$, $T^\Omega$ valides \\   
    \eSi{un des deux arbres n'est pas valide}{\Retour{\texttt{échec}}}{
        $(\mathbf{P}, \mathbf{Q}) \in [\mathbf{P}_{T^\Omega}] \times [\mathbf{Q}_{T^X}]$ \\
        \Retour{$\texttt{succès}, \mathbf{P},\mathbf{Q}, E_\textrm{bf}(\mathbf{A}, \mathbf{P}, \mathbf{Q})$}
    }
  \caption{Identifier $T^X$, $T^\Omega$ pour résoudre \eqref{eq:bfly-fact}.}
    \label{algo:heuristic}
\end{algorithm}

\smallparagraph{Résolution finale de \eqref{eq:bfly-fact}~}
\'Etant donné la matrice cible $\mathbf{A} := \tilde{\mathbf{Q}}^\top \tilde{\mathbf{A}} \tilde{\mathbf{P}}$, l'algorithme \ref{algo:heuristic} résout \emph{indépendamment} chaque problème \eqref{eq:min-partition} pour $\ell \in \set{L-1}$ via l'algorithme \ref{algo:partitionnement}. 
Si chaque problème \eqref{eq:min-partition} est bien résolu, et que les partitions trouvées forment des arbres $T^X$ et $T^\Omega$ valides dans le sens où ils satisfont les axiomes d'un \og cluster tree \fg, alors $\mathbf{A}$ possède la propriété de rang faible complémentaire pour $T^X$ et $T^\Omega$. 
On résout alors \eqref{eq:bfly-fact} via l'algorithme \eqref{algo:hierarchical}, en fixant $(\mathbf{P}, \mathbf{Q}) \in [\mathbf{P}_{T^\Omega}] \times [\mathbf{Q}_{T^X}]$.

\section{Expériences}
\label{sec:expérience}

Nous évaluons les performances empiriques de notre méthode pour factoriser %
$\mathbf{A} := \tilde{\mathbf{Q}}^\top \tilde{\mathbf{A}} \tilde{\mathbf{P}} + \epsilon (\| \tilde{\mathbf{A}} \|_F  / \| \mathbf{N} \|_F) \mathbf{N} $, où $\tilde{\mathbf{P}}$, $\tilde{\mathbf{Q}}$ sont des permutations aléatoires, $\mathbf{N}$ est une matrice avec des coefficients suivant une loi gaussienne centrée réduite, $\epsilon \geq 0$ %
contrôle le niveau relatif de bruit, et $\tilde{\mathbf{A}}$ correspond 
soit à une matrice butterfly orthogonale aléatoire définie par \cite{parker1995random}, soit à la matrice de la transformée de Fourier discrète (TFD). 
Nous appliquons l'algorithme \ref{algo:heuristic} avec $\{ \alpha_k \}_{k=1}^K := \{10^p \}_{p \in \{ -2, -1, 0, 1, 2\} }$ et $M=5$ sur 20 instances du problème pour $\epsilon \in \{ 0, 0.01, 0.03, 0.1 \}$, et pour $N \in \{ 2^L \}_{L \in \{ 2, \ldots, 7 \} }$. La taille $N > 128$ n'est pas considérée car l'algorithme \ref{algo:heuristic} a une complexité cubique en $N$: une exécution prend quelques minutes pour $N=64$, et une heure pour $N=128$.

\begin{figure}[t]
     \centering
    \includegraphics[width=\linewidth]{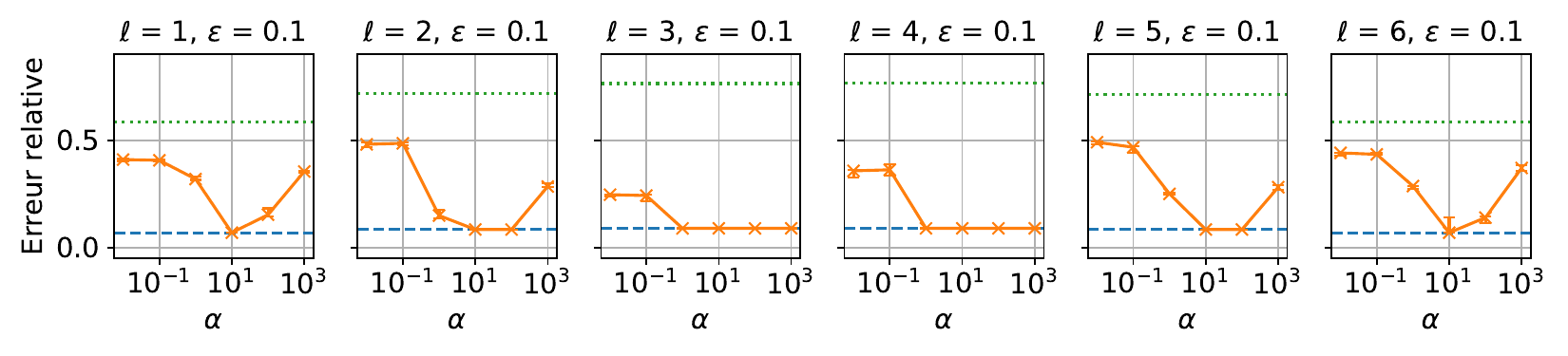}
    \caption{Résolution de \eqref{eq:min-partition} par l'algorithme \ref{algo:partitionnement} (avec $n=50$) pour une matrice butterfly orthogonale aléatoire de taille 128 (orange, plein). Les barres d'erreur (négligeables)
    affichent les extrema,  la croix indique la médiane. Bleu, tiret : erreur avec $T^X, T^\Omega$ connus. Vert, pointillé : erreur minimale sur 1000 tirages aléatoires de partitions, où l'on résout \eqref{eq:min-partition} en fixant les partitions.
    }
    \label{fig:reconstruction-partition}
    \vspace{-0.5em}
\end{figure}

\begin{figure}[t!]
    \centering
    \includegraphics[width=\linewidth]{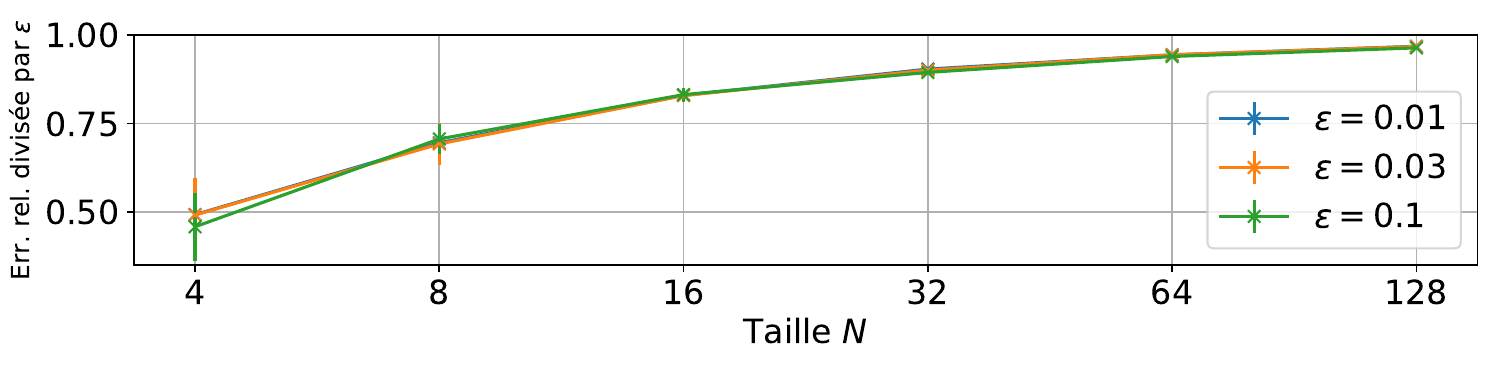}
    \caption{Erreur relative divisée par $\epsilon$ de l'algorithme \eqref{algo:heuristic} pour approcher une matrice butterfly orthogonale aléatoire bruitée (moyennes, écart-type).}
    \label{fig:error}
    \vspace{-1em}
\end{figure}

Lorsque $\tilde{\mathbf{A}}$ est une matrice butterfly orthogonale aléatoire, l'algorithme \ref{algo:heuristic} a 100 \% de succès sur les 20 instances du problème, pour tous les niveaux de bruit et toutes les tailles considérés, ce qui veut dire que l'algorithme \ref{algo:partitionnement} répété avec suffisamment de $\alpha$ et de graines aléatoires permet de résoudre indépendamment chaque problème \eqref{eq:min-partition} pour $\ell \in \set{L-1}$, et que les partitions trouvées forment des arbres $T^X$ et $T^\Omega$ valides. La figure \ref{fig:reconstruction-partition} illustre les résultats de l'algorithme \ref{algo:partitionnement}. En cas de succès, l'algorithme \ref{algo:heuristic} retourne la même erreur d'approximation que $\bflyerror{\mathbf{A}}{\tilde{\mathbf{P}}}{\tilde{\mathbf{Q}}}$ obtenue en connaissant $\tilde{\mathbf{P}}, \tilde{\mathbf{Q}}$. Dans le cas sans bruit, l'erreur relative atteint la précision machine. Dans le cas bruité, la figure \ref{fig:error} montre qu'elle est de l'ordre de $\epsilon$. 

Lorsque $\tilde{\mathbf{A}}$ est la matrice TFD, l'algorithme \ref{algo:heuristic} a également 100 \% de succès pour tous les niveaux de bruit et pour $N \leq 64$. Pour $N=128$, la taux de succès est à 100 \% pour le cas sans bruit, mais se dégrade dans le cas bruité comme illustré dans la table \ref{tab:taux-succès}. La robustesse pourrait être améliorée par exemple en résolvant les problèmes \eqref{eq:min-partition} pour $\ell \in \set{L-1}$ \emph{conjointement}, et non pas indépendamment.

\begin{table}[t]
    \caption{\label{tab:taux-succès}Taux de succès de l'algorithme \ref{algo:heuristic} sur 20 instances.}
    \vspace{-1em}
    \begin{center}
    \resizebox{\linewidth}{!}{
    \begin{tabular}{l*{5}{c}}
        \toprule
        $\epsilon$   &   0 &   0.01 &   0.03 &  0.1  \\
        \midrule
        Butterfly orthogonale aléatoire ($N=128$) &   100 \% &   100 \% &   100 \% &  100 \%  \\
        Transformée de Fourier discrète ($N=128$) &   100 \% &   95 \% &   90 \% &  50 \%  \\
        \bottomrule
    \end{tabular}
    }
    \vspace{-1.5em}
    \end{center}
\end{table}

\smallparagraph{Conclusion~}
Nous avons proposé une heuristique pour identifier sans hypothèse analytique les partitions d'une matrice en blocs de rang faible permettant une factorisation butterfly. 
Lever les verrous de robustesse et de passage à l'échelle de l'heuristique permettrait à terme de chercher une factorisation butterfly pour des opérateurs utilisés en traitement de signal ou en apprentissage, comme les transformées de Fourier sur graphes \cite{shuman2013emerging} ou les couches de réseaux de neurones \cite{dao2022monarch}.

\smallparagraph{Remerciements~} Ce travail a été soutenu par le projet ANR AllegroAssai ANR-19-CHIA-0009.

{\footnotesize
\begin{spacing}{0.8}
\bibliography{biblio}
\end{spacing}
}

\end{document}